\documentclass[11pt,leqno]{amsart}

\usepackage{enumerate}

\usepackage{mathrsfs}

\usepackage{amsthm}

\usepackage[english]{babel}

\usepackage{amsmath}

\usepackage{amssymb}

\usepackage{xypic}

\usepackage[all]{xy}

\usepackage{graphicx}
\usepackage{amsfonts}

\usepackage{enumerate}

\usepackage[T1]{fontenc}

\numberwithin{equation}{section}

\newcommand{\earrow}
{\mathrel{\text{\smash{\lower.477ex\hbox{$\overset{\text{\normalsize \smash{\lower.55ex\hbox{$\sim\,$}}}}{\to}$}}}}}

\newcommand{\bD}{\boldsymbol{\mathsf{D}}}

\newcommand{\rb}{\mathrm{b}}

\newcommand{\coh}{\mathrm{coh}}

\newcommand{\imin}[1]{#1^{-1}}

\newcommand{\op}{\mathrm{op}}

\newcommand{\sa}{\mathrm{sa}}

\newcommand{\rc}{{\R\textup{-c}}}

\newcommand{\XS}{X\times S}

\newcommand{\DXS}{\shd_{\XS/S}}

\newcommand{\lind}[1]{\underset{#1}{\varinjlim}}

\newcommand{\lpro}[1]{\underset{#1}{\varprojlim}}
\newcommand{\ot}{\sho^t}

\DeclareMathOperator{\codim}{codim}

\DeclareMathOperator{\rD}{\mathsf{D}}

\DeclareMathOperator{\id}{Id}\let\Id\id

\newcommand{\ds}{\displaystyle}

\numberwithin{equation}{section}

\DeclareMathOperator{\RH}{RH}

\newcommand{\C}{\mathcal{C}}

\newcommand{\D}{\mathcal{D}}

\newcommand{\T}{\mathcal{T}}

\renewcommand{\mod}{\mathrm{Mod}}

\newcommand{\psh}{\mathrm{Psh}}

\newcommand{\cov}{\mathrm{Cov}}

\newcommand{\CC}{\mathbb{C}}

\newcommand{\R}{\mathbb{R}}

\newcommand{\N}{\mathbb{N}}

\newcommand{\OO}{\mathcal{O}}

\newcommand{\I}{\mathrm{I}}

\newcommand{\osi}{\stackrel{\sim}{\gets}}

\newcommand{\iso}{\stackrel{\sim}{\to}}
\newcommand{\dbt}{\mathcal{D}\mathit{b}^t}

\newcommand{\db}{\mathcal{D}\mathit{b}}

\newcommand{\OW}{\OO^\mathrm{w}}

\newcommand{\CW}{\C^{{\infty ,\mathrm{w}}}}
\newcommand{\wtens}{\overset{\mathrm{w}}{\otimes}}

\newcommand{\rh}{\mathit{R}\mathcal{H}\mathit{om}}

\newcommand{\ho}{\mathcal{H}\mathit{om}}

\newcommand{\Ho}{\mathrm{Hom}}

\newcommand{\Rh}{\mathrm{RHom}}

\newcommand{\Lind}{\underrightarrow{\lim}}  %grazie Anna

\newcommand{\indl}[1]{\underset{#1}{``\underrightarrow{\mathrm{lim}}\mbox{''}}}

\newcommand{\exs}[3]{0 \to {#1} \to {#2} \to {#3} \to 0}

\newcommand{\towc}{\overset{wc}{\to}}

\def\shd{\mathcal{D}}

\def\shf{\mathcal{F}}\let\cF F
\let\cG G
\def\shh{\mathcal{H}}

\def\shm{\mathcal{M}}

\def\sho{\mathcal{O}}

\def\sht{\mathcal{T}}

\theoremstyle{plain}
\newtheorem{teo}{Theorem}[section]%[chapter]
\newtheorem{cor}[teo]{Corollary}%[chapter]
\newtheorem{prop}[teo]{Proposition}%[chapter]
\newtheorem{lem}[teo]{Lemma}%[chapter]
\theoremstyle{definition}
\newtheorem{es}[teo]{Example}%[chapter]
\newtheorem{oss}[teo]{Remark}%[chapter]
\newtheorem{df}[teo]{Definition}%[chapter]

\usepackage{latexsym}

\begin{document}

\title{Relative subanalytic sheaves II}

\date{}

\author{Teresa Monteiro Fernandes and Luca Prelli}

\thanks{The research of T. Monteiro Fernandes was supported by
Funda\c c{\~a}o para a Ci{\^e}ncia e Tecnologia, under the project: UIDB/04561/2020. The second author  is a member of the Gruppo Nazionale per l'Analisi Matematica, la Probabilit\`a e le loro Applicazioni (GNAMPA) of the Istituto Nazionale di Alta Matematica (INdAM).}



\maketitle
\vspace*{-\baselineskip}%

\begin{abstract}
We give a new construction of sheaves on a relative site associated to a product $X\times S$ where $S$ plays the role of a parameter space, expanding the previous construction by the same authors, where the subanalytic structure on $S$ was required. Here we let this last condition fall. In this way the construction becomes much easier to apply when dimension of $S$ is bigger than one. We also study the functorial properties of base change with respect to the parameter space.

\end{abstract}

\tableofcontents

\section{Introduction}

In a previous work (\cite{MFP14}), the authors introduced the notion of relative subanalytic sheaf associated to a sheaf on a subanalytic site of the form $(X\times S)_{sa}$, where $X$ and $S$ are real analytic manifolds. Here "relative" concerns the projection $p: X\times S\to S$ so that $S$ is regarded as a parameter space.

The background for all the constructions comes from the work of Kashiwara-Schapira \cite{KS5}.

The main purpose of \cite{MFP14} was to provide the subanalytic tools for a relative Riemann-Hilbert correspondence, now assuming that $X$ and $S$ are complex manifolds, generalizing the famous Kashiwara's Riemann-Hilbert correspondence in the absolute framework (the absolute case meaning that $S$ is a point). This correspondence was achieved in \cite{FMFS} under the assumption that $S$ is a complex curve. In that case the new objects are naturally endowed with an action of $p^{-1}\sho_S$-modules, $\shd_X$ being replaced by the sheaf of relative differential operators $\DXS$.

Our new goal is to provide the subanalytic tools to prove the relative Riemann-Hilbert correspondence for arbitrary dimension of $S$ (noted by $d_S$).
Let us explain the motivation:

While in the absolute case the reconstruction functor was based on the notion of "tempered distribution", "tempered holomorphic function", giving rise to the subanalytic sheaves $\db^t_X$ and $\sho_X^t$, in the relative case we needed to construct a "relative" version of those sheaves as sheaves on a relative site  "forgetting" the growth conditions on $S$. In \cite{MFP14}, our idea was to construct $X_{sa}\times S_{sa}$ respecting both subanalytic structures in $X$ and $S$, and there live the subanalytic sheaves $\db_{X\times S}^{\sharp,t,S}$ and $\sho_{X\times S}^{\sharp,t,S}$.

However it became clear that such a construction was too strict. Let us clarify the main obstruction:

Assume that $S$ is complex. We see that $p^{-1}\sho_S$, as a sheaf on $X_{sa} \times S_{sa}$ is not concentrated in degree zero unless $d_S=1$, since otherwise Stein subanalytic open subsets do not form a basis for $S_{sa}$. So, for arbitrary $d_S$, a new site should be constructed solving this difficulty.

The main idea was then to replace $X_{sa}\times S_{sa}$ by $X_{sa}\times S$ which means that we consider $S$ as a site with the usual topology. In this new setting we are able to consider infinite coverings of open subsets in $S$ and, in particular, we have a basis of the topology consisting of products of subanalytic open subsets of $X$ and Stein open subsets of $S$. Working with this construction is not trivial and requires some technical results (Propositions \ref{prop:imineta} and \ref{prop:cs2}).

In order to obtain them we introduce the notion of locally weakly quasi-compact site, which provides a general setting for many examples including the relative case.
Now the relative sheaves "live" on $X_{sa}\times S$ and we note $\db_{X\times S}^{t,S}$ the relative sheaf associated to $\db_{X\times S}^t$ and $\sho_{X\times S}^{t,S}$ the relative sheaf associated to $\sho_{X\times S}^{t}$.

We study the behavior of these sheaves under morphisms of the parameter space and we end with an application to a key particular case of the relative Riemann-Hilbert functor (the $S$-locally constant case studied in \cite{MFCS2}).

The contents of this paper are as follows.

In Section \ref{S:1} we prepare preliminaries to the study of the relative site. We recall some results about the subanalytic site $(X \times S)_{sa}$ and the site $X_{sa} \times S_{sa}$ to be used in the second part of the paper. After that we introduce the notion of locally weakly quasi compact site which will apply to the relative subanalytic site $X_{sa}\times S$. Using this generalization it is easier to establish some basic results about sheaves (in particular working with limits) and acyclic objects.

In Section \ref{S:2} we apply the previous framework to the case of a product $X\times S$ where $X$ and $S$ are complex manifolds. In the first part we re-state and re-prove the technical tools of \cite{MFP14} in this new framework. In order to do that we first need to prove some results to describe sections and acyclic objects in $X_{sa} \times S$ (Propositions \ref{prop:imineta} and \ref{prop:cs2}) and relate this site with $X \times S$ and $X_{sa} \times S_{sa}$ (Propositions \ref{prop:rhoS!} and \ref{prop:Tflabbycsoft}). After that we are able to study the relativization functor $(\bullet)^{RS}$ starting from the results of \cite{MFP14}. In the second part we study the functorial properties of the sheaf $\sho^{t,S}_{X\times S}$ in view of the Riemann-Hilbert correspondence.

\section{Preliminaries}\label{S:1}

In this section we state some results needed to study the relative site of Section \ref{S:2}. First we recall the constructions of \cite{KS5} and \cite{MFP14}. Then we define sheaves on locally weakly quasi-compact sites and we consider a useful family of acyclic objects. Throughout all the section, $k$ will denote a field.

\subsection{Sheaves on a subanalytic site}\label{subsection: Sheaves on a subanalytic site}

The results of this section  are extracted from \cite{KS5} (see also \cite{Pr08} for a more detailed study).

Let $X$ be a real analytic manifold and let $k$ be a field. Denote
by $\op(X_{sa})$ the category of open subanalytic subsets of $X$. One
endows $\op(X_{sa})$ with the following topology: $S \subset
\op(X_{sa})$ is a covering of $U \in \op(X_{sa})$ if for any
compact $K$ of $X$ there exists a finite subset $S_0\subset S$
such that $K \cap \bigcup_{V \in S_0}V=K \cap U$. We will call
$X_{sa}$ the subanalytic site.\\

Let $\mod(k_{X_{sa}})$ denote the category of sheaves on $X_{sa}$
and let $\mod_{\rc}(k_X)$ (resp. $\mod^c_{\rc}(k_X)$) be the abelian category of
$\R$-constructible sheaves on $X$ (resp. $\R$-constructible with compact support).\\

We denote by $\rho: X \to X_{sa}$ the natural morphism of sites.
We have functors

%\begin{eqnarray*}
%\rho_*:\mod(k_X) & \to & \mod(k_{X_{sa}}) \\
%\imin \rho:\mod(k_{X_{sa}}) & \to & \mod(k_X) \\
%\rho_!:\mod(k_X) & \to & \mod(k_{X_{sa}})
%\end{eqnarray*}

\begin{equation*}
\xymatrix{\mod_{\rc}(k_X) \subset \mod(k_X)   \ar@ <2pt>
[r]^{\hspace{1cm}\rho_*} &
  \mod(k_{X_{sa}}) \ar@ <2pt> [l]^{\hspace{1cm}\imin \rho}. }
\end{equation*}
The functors $\imin \rho$ and $\rho_*$ are the functors of inverse
image and direct image associated to $\rho$. The functor $\imin
\rho$ admits a left adjoint, denoted by $\rho_!$. The sheaf
$\rho_!F$ is the sheaf associated to the presheaf $\op(X_{sa}) \ni
U \mapsto F(\overline{U})$.
%In particular, for $U \in \op(X)$ one has $\rho_!k_U \simeq \lind {V
%\subset \subset U} \rho_*k_V$, where $V \in \op(X_{sa})$.

The functor $\rho_*$ is fully faithful and exact on
$\mod_{\rc}(k_X)$ and we identify $\mod_{\rc}(k_X)$ with its
image in $\mod(k_{X_{sa}})$ by $\rho_*$.\\

Thanks to $\R$-constructible sheaves we have some structure results for subanalytic sheaves, namely

\begin{itemize}
\item[(i)] Let $G \in \mod^c_{\rc}(k_X)$ and let
$\{F_i\}$ be a filtrant inductive system in $\mod(k_{X_{sa}})$.
Then we have the isomorphism
$$\lind i \,\Ho_{k_X}(\rho_*G,F_i) \iso
\Ho_{k_{X_{sa}}}(\rho_*G,\lind i\, F_i).$$
In particular, if $U \in \op(X_{sa})$ is relatively compact
$$\lind i \,\Gamma(U;F_i) \iso
\Gamma(U;\lind i\, F_i).$$

\item[(ii)] Let $F \in \mod(k_{X_{sa}})$. There exists a small filtrant
inductive system $\{F_i\}_{i \in I}$ in $\mod_{\rc}^c(k_X)$ such
that $F \simeq \lind i\, \rho_*F_i$.
\end{itemize}

In order to find suitable acyclic resolutions, let us recall the definition of quasi-injective sheaves. $F \in \mod(k_{X_{sa}})$ is quasi-injective, if for any relatively compact $U \in \op(X_{sa})$ the restriction morphism $\Gamma(X;F) \to \Gamma(U;F)$ is surjective.

Let  $G \in \mod_{\rc}(k_X)$. Then the following hold:
\begin{itemize}
\item[(i)]
 The family of quasi-injective sheaves is injective with respect to
the functor $\Ho(G,\bullet )$. In particular, it is injective with respect to the functor $\Gamma(U;\bullet)$, for each $U \in \op(X_{sa})$.

\item[(ii)]
The family of quasi-injective sheaves is injective with
respect to the functor $\ho(G,\bullet )$. In particular, it is injective with respect to the functor $\Gamma_U(\bullet)$, for each $U \in \op(X_{sa})$.

\end{itemize}

\subsection{Sheaves on a product of subanalytic sites}\label{subsection: Sheaves on a product of subanalytic sites}

 As a motivation for the contents of this section we start by recalling the example provided by the construction of \cite{MFP14}. Let $X$ and $S$ be two real analytic manifolds. It is possible to consider the product of sites $X_{sa} \times S_{sa}$. A basis $\T$ of this topology consists of products $U \times V$ where $U \in \op(X_{sa})$ and $V \in \op(S_{sa})$ are relatively compact. Then $W \in \op(X_{sa} \times S_{sa})$ if it is a locally finite union of elements of $\T$. A subset $T \subset \op(X_{sa} \times S_{sa})$ is a covering of $U \times V \in \T$ if it admits a finite subcover. More generally, $T \in \cov(W)$, where $W \in \op(X_{sa} \times S_{sa})$, if and only, if for any $\T \ni U \times V \subset W$, one has $\{W' \cap (U \times V),\,\ W' \in T\} \in \cov(U \times V)$ (i.e. the intersection of $T$ with $U \times V$ admits a finite subcover).
We denote by  $\rho': X \times S \to X_{sa} \times S_{sa}$  the natural morphism of sites.
There are well defined functors:
$$
\xymatrix{\mod(k_{X \times S})
\ar@ <2pt> [r]^{\mspace{0mu}\rho'_*} &
  \mod(k_{X_{sa} \times S_{sa}}) \ar@ <2pt> [l]^{\mspace{0mu}\rho'{}^{-1}}. }
$$
The functor $\rho'_*$ is fully faithful. The functor $\rho'{}^{-1}$ admits  a left adjoint, denoted by $\rho'_!$. It is fully faithful, exact and commutes with $\otimes$. Given $F \in \mod(k_{X_{sa} \times S_{sa}})$, $\rho'_!F$ is the sheaf associated to the presheaf $$\T \ni U \times V \mapsto F(\overline{U} \times \overline{V})$$

Let us consider the category $\coh(\T)$ of sheaves of
$\mod(k_{X \times S})$, admitting a finite resolution by sums $\oplus_{i \in I} k_{U_i \times V_i}$ with $I$
finite and $U_i \times V_i\in \T$ for each $i$. The category  $\coh(\T)$ is additive and stable by kernels and cokernels.

Moreover:
\begin{itemize}
\item{The restriction of $\rho'_*$ to $\coh(\T)$ is exact. Hence the latter can be regarded as a subcategory of $\mod(k_{X_{sa} \times S_{sa}})$.}
\item{A sheaf $F \in \mod(k_{X_{sa} \times S_{sa}})$ can be seen as a filtrant inductive limit $\lind i\, \rho'_*F_i$ with $F_i \in \coh(\T)$.}
\item{The functors $\Ho(G,\bullet)$, $\ho(G,\bullet)$, with $G \in \coh(\T)$, commute with filtrant $\Lind$.}
\end{itemize}
Finally, recall (cf \cite{EP}) that $F \in \mod(k_{X_{sa} \times S_{sa}})$ is $\T$-flabby if the restriction morphism $\Gamma(X;F) \to \Gamma(W;F)$ is surjective for each $W \in \T$. $\T$-flabby objects are $\Ho(G,\bullet)$-acyclic and $\ho(G,\bullet)$-acyclic for each $G \in \coh(\T)$.

\subsection{Sheaves on locally weakly quasi-compact sites}\label{subsection: Sheaves on locally weakly quasi-compact sites}

Let $\C_X \subset \op(X)$, where $\op(X)$ is the category open subsets of a topological space $X$. Until the end of the section we will
assume that $\C_X$ is stable under finite unions and intersections. As a category, if $U,V \in \C_X$ one has $U \to V$ if and only if $U \subset V$.
If $U \in \C_X$ we denote by $\C_U$ the subcategory of $\C_X$ consisting of open subsets $V \subset U$.
We suppose that a Grothendieck topology on $\C_X$ is given and for each $U \in \C_X$ the set $\cov(U)$ is given by a subfamily of topological coverings of $U$. We denote by $X$ the associated site. Given a covering $S$ of $U$, and $V \to U$, we set $S \times_U V = \{W \cap V,\,\ W \in S\}$.

\begin{df}\label{df:Xf} One defines the site $X^f$ as follows: $S$ is a covering
of $U \in \C_X$ in $X^f$ if it admits a finite refinement. We set $\cov^f(U)$ the family of coverings of $U$ in $X^f$.
\end{df}

\begin{df}\label{df:wqc} (i) Let $U,V \in \C_X$ and let $V \to U$. One says that $V$  is weakly
  quasi-compact (= wc) in $U$ if, for any  $S \in \cov(U)$, we
  have $S \times_U V \in \cov^f(V)$.
We will write $V \towc U$ to say that $V \to U$ is  wc in $U$.

(ii) Let $U,V \in \C_X$ with $V \to U$ and let $T\in \cov(V)$, $S \in \cov(U)$.  We write $T
\towc S$ if $T$ is a refinement of $S \times_U V$ and for any $V' \in T$, $U' \in S$ one has $V' \to U'$ iff $V'\towc U'$.
\end{df}

Note that when $S=U \in \C_U$ and $T=V \in \C_V$ we recover
Definition \ref{df:wqc} (i).

\begin{df}\label{df:lwqc}
A site $X$ is locally weakly quasi-compact if:
\begin{itemize}
\item[LWQC1] for each $U \in \C_X$ $\{V \towc U\} \in
\cov(U)$,
\item[LWQC2] for each $U,V,W \in \C_X$, if $W \towc U$ and $W \towc V$ one has $W \towc U \times_X V$,
\item[LWQC3] for each $V \towc U$ there exists $W \towc U$ with $V \towc W$.
\end{itemize}
\end{df}

\begin{lem}\label{lem:TwcS} Let $V \towc U$. Then for each $S \in \cov(U)$ there exists $T^f \in \cov^f(V)$ such that $T^f \towc S$.
\end{lem}
\begin{proof} By LWQC3 there exists $V' \towc U$ with $V \towc V'$. By LWQC1 for each $U_j \in S$ we have $\{W_{ij} \towc U_j \times_U V'\} \in \cov(U_j \times_U V')$, hence $\bigcup_j\{W_{ij} \towc U_j \times_U V'\} \in \cov(V')$ is a refinement of $S \times_U V'$. We have $T:=\bigcup_j\{W_{ij} \towc U_j \times _UV' \,\to V'\} \in \cov(V')$ and $T^f:=T \times_{V'} V \in \cov^f(V)$. This implies $T^f \towc S$.
\end{proof}

\begin{lem}\label{lem:fac+} Let $F \in \psh(k_X)$, and let $U \in \C_X$. If $F$ is a sheaf on $X^f$, then  for any $V \towc U$ the morphism
\begin{equation}\label{eq:factors+}
F^+(U) \to F^+(V)
\end{equation}
 factors through $F(V)$.
\end{lem}
\begin{proof}
Let $S \in \cov(U)$. There is a finite refinement $T^f \in
\cov^f(V)$ of $S \times_U V$. Then the morphism \eqref{eq:factors+} is
defined by
\begin{eqnarray*}
F^+(U) & \simeq & \lind {S \in \cov(U)} \hspace{-3mm} F(S)\\
 & \to & \lind {S \in \cov(U)} \hspace{-3mm}F(S \times_U V )\\
& \to & \lind {T^f \in \cov^f(V)}\hspace{-3mm} F(T^f)\\
& \to & \lind {T \in \cov(V)}\hspace{-3mm}F(T)\\
& \simeq & F^+(V).
\end{eqnarray*}
$F$ is a sheaf on $X^f$, hence $F(T^f) \simeq F(V)$ and the result follows.
\end{proof}

\begin{cor}\label{cor:fac|} Under the hypothesis of Lemma \ref{lem:fac+}, let $S \in \cov(U)$ and $T \in
\cov(V)$. If $T \towc S$, then the morphism
\begin{equation}\label{eq:factorscov+}
F^+(S) \to F^+(T) \end{equation}
 factors through $F(T)$. In particular, if $T \in \cov^f(V)$, then the
 morphism \eqref{eq:factorscov+} factors through $F(V)$.
\end{cor}

From now on we assume that $X$ is locally weakly quasi-compact.

\begin{lem}\label{lem:lproc} Let $F \in \psh(k_X)$. Then we have the isomorphism
$$
\lpro {V\towc U} \lind {V\towc W\towc U} F(W) \iso \lpro {V \towc U} F(V).
$$
\end{lem}
\begin{proof} The result follows since, for any $U\in\C_X$, for each $V \towc U$, there exists $W \towc U$ such that $V \towc W$ by LWQC3. Let $V\towc U$. The restriction morphism $F(U) \to F(V)$ factors through $\ds{\lind {W \towc V \towc U}F(W)}$. Taking the projective limit we obtain the result.
\end{proof}

\begin{lem}\label{lem:fac++} Let $F \in \psh(k_X)$, and let $U \in \C_X$. If $F$ is a sheaf on $X^f$, then  for any $V \towc U$ the morphism
\begin{equation}\label{eq:factors++}
F^{++}(U) \to F^{++}(V)
\end{equation}
 factors through $F(V)$.
\end{lem}
\begin{proof}
Let $S \in \cov(U)$. Thanks to Lemma \ref{lem:TwcS} we can construct a refinement $T^f \in \cov^f(V)$ of $S \times_U V$ such that $T^f \towc S$.
The morphism \eqref{eq:factors++} is
defined by
\begin{eqnarray*}
F^{++}(U) & \simeq & \lind {S \in \cov(U)} \hspace{-3mm} F^+(S)\\
 & \to & \lind {S \in \cov(U)} \hspace{-3mm}F^+(S \times_U V )\\
& \to & \lind {T^f \in \cov^f(V)}\hspace{-3mm} F^+(T^f)\\
& \to & \lind {T \in \cov(V)}\hspace{-3mm}F^+(T)\\
& \simeq & F^{++}(V).
\end{eqnarray*}
By Corollary \ref{cor:fac|}, the morphism
$$
\lind {S \in \cov(U)} \hspace{-3mm} F^+(S) \to  \lind {T^f \in \cov^f(V)}\hspace{-3mm} F^+(T^f)
$$
 factors through $F(V)$ and the result follows.
\end{proof}

\begin{cor}\label{cor:faclimlim} Let  $F \in \psh(k_X)$. If $F$ is a sheaf on $X^f$, then:
\begin{itemize}
\item[(i)] for any $U \in \C_X$ one has the isomorphism $\ds{\lind {V \towc U}F(U)
\iso \lind {V \towc U}F^{++}(U)}$. \\
\item[(ii)] for any $U \in \C_X$ one has the isomorphism $\lpro {V \towc U}F(V)
\iso \lpro {V \towc U}F^{++}(V)$.
\end{itemize}
\end{cor}
\begin{proof}
 By Lemma \ref{lem:fac++} for each $V
\towc U$ we have a commutative diagram
$$
\xymatrix{ F^{++}(U) \ar[r]  \ar[rd] & F^{++}(V)  \\
F(U) \ar[u] \ar[r] & F(V) \ar[u] .}
$$
 Passing to the inductive (resp. projective) limit we obtain $\lind {V \towc U}F(U) \iso \lind {V \towc U}F^{++}(U)$ (resp. $\lpro {V \towc U}F(V)
\iso \lpro {V \towc U}F^{++}(V)$) as required. The isomorphism follows from the universal property and the uniqueness, up to isomorphisms, of the injective (resp. projective limit).
\end{proof}

Let $\{F_i\}_{i \in I}$ be a filtrant inductive system in
$\mod(k_X)$. One sets
\begin{equation*}
\begin{array}{l}
\text{$\indl i\, F_i$ = inductive limit in the category of presheaves,} \\
\text{$\lind i\, F_i$ = inductive limit in the category of sheaves.}
\end{array}
\end{equation*}
Recall that $\lind i\ F_i=(\indl i F_i)^{++}$.

\begin{prop}\label{prop:limXf} Let $\{F_i\}_{i \in I}$ be a filtrant inductive
system in $\mod(k_X)$. Then $\indl i
F_i$ is a sheaf on $X^f$.
\end{prop}
\begin{proof} Let $U
\in \C_X$ and let $S$ be a finite covering of $U$. Since
$\lind i$ commutes with finite projective limits we obtain the
isomorphism $(\indl i F_i)(S) \iso \lind i \,F_i(S)$ and $F_i(U)
\iso F_i(S)$ since $F_i \in \mod(k_\T)$ for each $i$.
Moreover the family of finite coverings of $U$ is cofinal in
$\cov(U)$. Hence $\indl i\, F_i \iso (\indl i F_i)^+$ on $X^f$. Applying once
again the functor $(\bullet)^+$ we get
$$\indl i F_i \simeq (\indl i F_i)^+ \simeq (\indl i F_i)^{++} \simeq \lind i \,F_i$$
on $X^f$ and the result follows.
 \end{proof}

Now, setting $F=\indl i F_i$ in Lemma \ref{lem:fac++} and in Corollary \ref{cor:faclimlim} we obtain the following results:

\begin{prop}\label{prop:fac} Let $\{F_i\}_{i \in \I}$ be a filtrant inductive system in $\mod(k_X)$ and let $U \in \C_X$.
Then for any $V \towc U$ the morphism
\begin{equation*}
\Gamma(U;\lind i\, F_i) \to \Gamma(V;\lind i\, F_i)
\end{equation*}
 factors through $\lind i \Gamma(V;F_i)$.
\end{prop}

\begin{cor}\label{cor:faclim} Let $\{F_i\}_{i \in \I}$ be a filtrant inductive system in $\mod(k_X)$.
\begin{itemize}
\item[(i)] For any $U \in \C_X$ one has the isomorphism $$\lind {U
\towc V,i}\Gamma(V;F_i) \iso \lind {U \towc
V}\Gamma(V; \lind i\, F_i)$$ \\
\item[(ii)] For any $U \in \C_X$ one has the isomorphism $$\lpro {V \towc
U}\lind i\, \Gamma(V;F_i) \iso \lpro {V \towc
U}\Gamma(V;\lind i \,F_i)$$
\end{itemize}
\end{cor}

%\begin{es} Examples of locally weakly compact sites are:
%\begin{enumerate}
%\item Locally compact topological spaces, relatively compact open subsets are locally weakly compact.
%\item Locally compact spaces with an action of $\mu$ of $\RP$ endowed with the conic topology (i.e. open subsets are $\mu$-invariant), orbits $\mu(U,\RP)$ of relatively compact open subsets are weakly compact.
%\item Noetherian spaces, quasi-compact open subsets are locally weakly compact.
%\item The subanalytic site of \cite{KS5}, relatively compact subanalytic open subsets are weakly compact.
%\item The conic subanalytic site of \cite{Pr11} (see also \cite{Pr13}), orbits $\mu(U,\RP)$ of relatively compact subanalytic open subsets are weakly compact.
%\item More generally, the $\T$-sites of \cite{KS5} and the $\T_{loc}$ sites of \cite{EP}.
%\item The $\T$-spectrum of \cite{EP}, quasi-compact open subsets are weakly compact.
%\item The relative subanalytic site of \cite{MFP14}, products of relatively compact subanalytic open subsets are weakly compact.
%\item The new relative site $X_{sa} \times S$ of Section \ref{S:2}, products of relatively compact subanalytic open subsets of $X$ and relatively compact open subsets of $S$ are weakly compact.
%\end{enumerate}
%\end{es}

\subsection{c-soft sheaves}

Let $X$ be a locally weakly quasi-compact site.

\begin{df} We say that a sheaf $F$ on $X$ is c-soft if the restriction
morphism $\Gamma(W;F) \to \lind{V \towc U}\Gamma(U;F)$ is
surjective for each $V,W \in \C_{Y}$ ($Y \in \C_X$) with $V,W \towc Y$ and $V\towc W$.
\end{df}

\begin{prop}\label{prop:cs1} Let $\exs{F'}{F}{F''}$ be an exact sequence in
$\mod(k_X)$, and assume that $F'$ is c-soft. Then the sequence
$$\exs{\lind{V \towc U}\Gamma(U;F')}{\lind{V \towc U}\Gamma(U;F)}{\lind{V \towc U}\Gamma(U;F'')}$$
is exact for any $V \in \C_X$.
\end{prop}
\begin{proof} Let $s'' \in \lind{V \towc U}\Gamma(U;F'')$. Then
there exists $V \towc U$ such that $s''$ is represented
by $s''_U \in \Gamma(U;F'')$. Let $\{U_i\}_{i \in I} \in \cov(U)$
such that there exists $s_i \in \Gamma(U_i;F)$ whose image is
$s''_U|_{U_i}$ for each $i$. There exists $W \towc U$ with
$V\towc W$, a finite covering $\{W_j\}_{j=1}^n$ of $W$
and a map $\varepsilon:J \to I$ of the index sets such that $W_j
\towc U_{\varepsilon(j)}$. We may argue by induction on
$n$. If $n=2$, set $U_i=U_{\varepsilon(i)}$, $i=1,2$. Then
 $(s_1-s_2)|_{U_1 \cap U_2}$ belongs to $\Gamma(U_1 \cap U_2;F')$,
 and  its restriction defines an element of $\lind {W_1 \cap
 W_2\towc W'}\Gamma(W';F')$, hence it extends to $s' \in \Gamma(U;F')$. By
 replacing $s_1$ with $s_1-s'$ on $W_1$ we may assume that
 $s_1=s_2$ on $W_1 \cap W_2$. Then there exists $s \in \Gamma(W_1
 \cup W_2;F)$ with $s|_{W_i}=s_i$. Thus the induction proceeds.
 \end{proof}

\begin{prop}\label{prop:www} Let $\exs{F'}{F}{F''}$ be an exact sequence in
 $\mod(k_{X})$, and assume $F',F$ c-soft. Then $F''$
 is c-soft.
\end{prop}
\begin{proof} Let $V,W,Y \in \C_X$ with $V\towc W \towc Y$ and let us
consider the diagram below
$$ \xymatrix{\Gamma(W;F) \ar[d]^\alpha \ar[r] & \Gamma(W;F'') \ar[d]^\gamma \\
\lind {V \towc U}\Gamma(U;F) \ar[r]^\beta & \lind {V \towc U}\Gamma(U;F'').} $$  The morphism
$\alpha$ is surjective since $F$ is c-soft and $\beta$ is
surjective by Proposition \ref{prop:cs1}. Then $\gamma$ is
surjective.
\end{proof}

\begin{prop} Let $V \in \C_X$. The family of c-soft sheaves is injective with respect
to the functor $\lind {V \towc U}\Gamma(U;\bullet)$ .
\end{prop}
\begin{proof} The family of c-soft sheaves contains injective sheaves, hence it is cogenerating. Then the result follows from Propositions \ref{prop:cs1} and \ref{prop:www}.
\end{proof}

\begin{oss} Remark that when $X$ is a locally compact space and we consider the category $\mod(k_X)$, the above definition of c-soft sheaves coincides with the classical one.
\end{oss}

%\begin{es} Examples of c-soft sheaves:
%\begin{enumerate}
%\item When $X$ is a locally compact topological space: $F$ is c-soft if $F(X) \to F(K)$ is surjective for each $K \subset X$ compact.
%\item When $X$ is a locally compact topological space endowed with the conic topology: $F$ is c-soft if $F(X) \to \lind {\mu(U,\RP) \subset\subset V}F(V)$ is surjective for each $U$ relatively compact open subset of $X$ (with the natural topology of $X$) and $V$ is an open cone.
%\item When $X$ is a Noetherian space: $F$ is c-soft if $F(X) \to F(K)$ is surjective for each $U \subset X$ open quasi-compact.
%\item When $X$ is the semialgebraic site of \cite{De}, the subanalytic site of \cite{KS5}, the conic subanalytic site of \cite{Pr11}, the relative site of \cite{MFP14} or, in general, the $\T$-site of \cite{EP}, c-soft sheaves and $\T$-flabby sheaves coincide.

%\item The new relative site of Section \ref{S:2}:  $F$ is c-soft if $F(U) \to F(K)$ is surjective for each $K \subset X$, where $K$ is a finite union of products of definably compact subsets of $S$ and relatively compact open subanalytic subsets of $X$ and $U$ is relatively compact with $K \subset U$.
%\end{enumerate}
%\end{es}

\section{Relative subanalytic site}\label{S:2}

This section resumes the purpose of \cite{MFP14}, allowing the dimension $d_S$ of the parameter space to be arbitrary. From now on we set $k=\CC$.

\subsection{Relative subanalytic sheaves}

Let $X$ and $S$ be two real analytic manifolds. Let us denote by $p:X\times S\to S$ the projection. We consider the following Grothendieck topologies on $X \times S$.
\begin{itemize}
\item The one generated by the usual topology.
\item The subanalytic topology of \cite{KS5} and \cite{Pr08}, defining the site $(X \times S)_{sa}$, where $\op((X \times S)_{sa})$ are open subanalytic subsets of $X \times S$ and the coverings are those admitting a locally finite refinement. See Section \ref{subsection: Sheaves on a subanalytic site} for more details.
\item The product of sites $X_{sa} \times S_{sa}$ of \cite{MFP14}. This is the site where the open sets are locally finite unions of elements of the family $\T$ of products $U\times V$ with $U$ and $V$ relatively compact subanalytic open respectively in $X$ and $S$, and the coverings are the families of such open sets admitting a locally finite refinement. See Section \ref{subsection: Sheaves on a product of subanalytic sites} for more details.
\item The product of sites $X_{sa} \times S$ for which the elements of $\op(X_{sa} \times S)$ are finite unions of products $U \times V$, with $U \in \op(X_{sa})$ and $V \in \op(S)$. With this topology, $X_{sa} \times S$ is a locally weakly quasi-compact site. See Section \ref{subsection: Sheaves on locally weakly quasi-compact sites} for more details.
\end{itemize}
Let us describe the coverings in $X_{sa} \times S$:
  %Let $p_X: X \times S \to X$ and $p_S: X \times S \to S$ be the projections.
$T \subset \op(X_{sa} \times S)$ is a covering of $W=U \times V \in \op(X_{sa} \times S)$ if and only if it admits a refinement $\{U_i \times V_j\}_{i \in I, j \in J}$ with $\{U_i\}_{i \in\ I} \in \cov(U)$ (in $X_{sa}$) and $\{V_j\}_{j \in J} \in \cov(V)$ (in $S$). In particular, when $U$ is relatively compact $I$ is finite. %If $W=\cup_{j \in J}W_j$ with $W_j=U_j \times V_j \in \op(X_{sa} \times S)$, then $T$ is a covering of $W$ if $T \times W_j$ is a covering of $W_j$ for each $j$. It is a locally weakly compact site and $U \times V \towc X \times S$ if $U \subset\subset X$ and $V \subset\subset S$.
For a general $W \in \op(X_{sa} \times S)$, $T \subset \op(X_{sa} \times S)$ is a covering of $W$ if for each $U \times V \in \op(X_{sa} \times S)$ with $U \times V \subset W$ one has $T \times_W (U \times V) \in \cov(U \times V)$.

Let us make explicit the notion of weakly quasi-compact in $X_{sa} \times S$: let $W, W' \in \op(X_{sa} \times S)$ and suppose that $W=\cup_{j \in J}(U_j \times V_j)$ with $J$ finite. Then $W \towc W'$ if $U_j$ and $V_j$ are relatively compact and $U_j \times \overline{V}_j \subset W'$ for each $j \in J$.

We have the following commutative diagram, where the arrows are natural morphisms of sites induced by the inclusion of families of open subsets.
\begin{equation}\label{eq:diagram sites}
\xymatrix{
 & X_{sa} \times S \ar[dr]^a &  \\
 X \times S \ar[rr]^{\rho'} \ar[ur]^{\rho_S} \ar[dr]^\rho & & X_{sa} \times S_{sa} \\
 & (X \times S)_{sa} \ar[ur]^\eta &
}
\end{equation}

It follows from Theorem 3.9.2 of \cite{Ta94} that the functor $a_*$ is fully faithful and $\imin a \circ a_* \iso \id$. Similarly, the functor $\rho_{S*}$ is fully faithful and $\imin{\rho_S} \circ \rho_{S*} \iso \id$.

\begin{prop}\label{prop:imineta} Let $U \times V \in \op(X_{sa} \times S)$ and let $F \in \mod(\CC_{X_{sa} \times S_{sa}})$. Then
$$
\Gamma(U \times V;\imin aF) \simeq \lpro {V' \subset\subset V}\Gamma(U \times V';F)
$$
with $V'$ subanalytic.
\end{prop}
\begin{proof}
First suppose that $U$ is relatively compact. Remark that $\{V' \subset\subset V, V' \in \op(S_{sa})\}$ is a covering of $V$ and
$$
\Gamma(U \times V;\imin aF) \simeq \lpro {V' \subset\subset V}\Gamma(U \times V';\imin aF)
$$
with $V'$ subanalytic. There exists a filtrant family $\{F_i\}$ of $\T$-coherent sheaves such that $F \simeq \lind i\, \rho'_*F_i$. Remark that $\rho' = a \circ \rho_S$ and we have
$$
\imin a\lind i\,\rho'_*F_i \simeq \lind i \,\imin a \circ a_* \circ \rho_{S*}F_i \simeq \lind i \,\rho_{S*}F_i.
$$
We have the chain of isomorphisms
\begin{eqnarray*}
\Gamma(U \times V;\imin aF) & \simeq & \lpro {V' \subset\subset V}\Gamma(U \times V';\imin aF) \\
& \simeq & \lpro {V' \subset\subset V}\Gamma(U \times V';\imin a\lind i\,\rho'_*F_i) \\
& \simeq & \lpro {V' \subset\subset V}\Gamma(U \times V';\lind i\, \rho_{S*}F_i) \\
& \simeq & \lpro {V' \subset\subset V}\lind i\,\Gamma(U \times V';\rho_{S*}F_i) \\
& \simeq & \lpro {V' \subset\subset V}\lind i\,\Gamma(U \times V';\rho'_*F_i) \\
& \simeq & \lpro {V' \subset\subset V}\Gamma(U \times V';\lind i\,\rho'_*F_i) \\
& \simeq & \lpro {V' \subset\subset V}\Gamma(U \times V';F)
\end{eqnarray*}
with $V'$ subanalytic. The fourth and sixth isomorphisms follow from Corollary \ref{cor:faclim} and the fifth follows since $U$ and $V'$ are subanalytic.

%First let us suppose that $U$ is relatively compact. Let $V' \subset\subset V'' \subset\subset V$ with $V',V''$ open subanalytic. By Lemma \ref{lem:fac++} the restriction $\Gamma(U \times V'';\imin aF) \to \Gamma(U \times V';\imin aF)$ factors through $\Gamma(U \times V';F)$. Then we have the isomorphism
%$$
%\Gamma(U \times V;\imin aF) \simeq \lpro {V' \subset\subset V}\Gamma(U \times V';F)
%$$
Suppose now that $U$ is not relatively compact.  Let $\{U_n\}_{n \in \N}$ be a covering of $X$ consisting of open subanalytic subsets such that $U_n \subset\subset U_{n+1}$ for every $n \in \N$. Then
\begin{eqnarray*}
\Gamma(U \times V;\imin aF) & \simeq & \lpro n\,\Gamma(U \cap U_n \times V; \imin aF)\\
& \simeq & \lpro {n,V' \subset\subset V}\Gamma(U \cap U_n \times V';F) \\
& \simeq & \lpro {V' \subset\subset V}\Gamma(U \times V';F).
\end{eqnarray*}
\end{proof}

\begin{cor}\label{cor:imineta} Let $W \in \op(X_{sa} \times S)$ and let $F \in \mod(\CC_{X_{sa} \times S_{sa}})$. Then
$$
\Gamma(W;\imin aF) \simeq \lpro {W' \towc W}\Gamma(W';F)
$$
with $W' \in \op(X_{sa} \times S_{sa})$.
\end{cor}

%\begin{cor} The adjunction morphism induces the isomorphism $\imin a \circ a_* \iso \id$.
%\end{cor}

We shall need to introduce a left adjoint to the functor $\imin {\rho_S}$.

\begin{prop} \label{prop:rhoS!} The functor $\imin {\rho_S}$ has a left adjoint, denoted by $\rho_{S!}$. The functor $\rho_{S!}$ is exact and commutes with tensor products.
\end{prop}
\begin{proof}
With the notations of \eqref{eq:diagram sites} $\rho'{}^{-1} \circ a_* \simeq \imin {\rho_S} \circ \imin a \circ a_* \simeq \imin {\rho_S}$, where the second isomorphism follows from Theorem 3.9.2 of \cite{Ta94}. By Proposition 2.4.3 of \cite{EP} the functor $\rho'{}^{-1}$ has a left adjoint $\rho'_!$, so the composition $\rho'{}^{-1} \circ a_* \simeq \imin {\rho_S}$ admits a left adjoint $\imin a \circ \rho'_! =: \rho_{S!}$. Remark that by definition (see also Proposition 2.4.4 of \cite{EP}) $\rho_{S!}$ is exact and commutes with tensor products.
\end{proof}

\begin{oss} Thanks to Proposition \ref{prop:imineta} and the fact that $\imin a \circ \rho'_! =: \rho_{S!}$, one can prove that $\rho_{S!}F$ (with $F \in \mod(\CC_X)$) is the sheaf associated to the presheaf $U \times V \mapsto \lind {U'}\,\Gamma(U' \times V;F)$, with $ {\overline{U} \subset U'}$. An alternative proof can be obtained adapting the proof of Proposition 6.6.3 of \cite{KS5}.
\end{oss}

\begin{prop}\label{prop:cs2} The category of c-soft sheaves is injective with respect
to the functor $\Gamma(U \times V; \bullet)$, $U \in \op(X_{sa})$, $V \in \op(S)$.
\end{prop}
\begin{proof} First remember that, if $U' \subset\subset X$ and $V' \subset\subset V$ ($U',V'$ subanalytic), then $U' \times V'$ is weakly compact in $U \times V$ in the relative site $X_{sa} \times S$.  Let $\{U_n\}_{n \in \N}$ (resp. $\{V_n\}_{n \in\N}$) be a covering of $X$ (resp. $V$) consisting of open subanalytic subsets such that $U_n \subset\subset U_{n+1}$ (resp. $V_n \subset\subset V_{n+1})$ for every $n \in \N$. For any $G \in \mod(\CC_{X_{sa} \times S})$ we have
$$
\Gamma(U \times V;G) \simeq \lpro n\, \Gamma(U \cap U_n \times V_n;G) \simeq \lpro n\, \lind {V'_n \supset\supset V_n} \Gamma(U \cap U_n \times V'_n;G).
$$
Set for short $W_n:=U \cap U_n \times V'_n$, $n \in \N$ and
$$\lind {W_n}\,\Gamma(W_n;G):=\lind {V'_n \supset\supset V_n} \Gamma(U \cap U_n \times V'_n;G).$$
Then $U \cap U_n \times V_n \towc W_n$, $W_n$ is weakly compact and we can use the results of Section \ref{subsection: Sheaves on locally weakly quasi-compact sites}. Remark that $U \cap U_n$ is relatively compact, hence every covering in $X_{sa}$ admits a finite refinement.

 Take an exact sequence $\exs{F'}{F}{F''}$, and suppose
$F'$ c-soft.  All the
sequences
$$
\exs{\lind {W_n}\, \Gamma(W_n;F')}{\lind {W_n}\, \Gamma(W_n;F)}{\lind {W_n}\, \Gamma(W_n;F'')}$$
are exact by Proposition \ref{prop:cs1}, and the morphism
$$\lind {W_{n+1}}\, \Gamma(W_{n+1};F') \to
\lind {W_n}\, \Gamma(W_n;F')$$
is surjective for all
$n$. Then by Proposition 1.12.3 of \cite{KS1} the sequence
$$\exs{\lpro n\,\lind {W_n}\, \Gamma(W_n;F')}{\lpro n\,\lind {W_n}\, \Gamma(W_n;F)}{\lpro n\,\lind {W_n}\, \Gamma(W_n;F'')}$$ is exact.
\end{proof}

\begin{prop}\label{prop:Tflabbycsoft} Let $F \in \mod(\CC_{X_{sa} \times S_{sa}})$ be $\T$-flabby. Then $\imin aF$ is c-soft.
\end{prop}
\begin{proof} Let us first prove that, if $W \towc W'$ with $W, W'\in\op(X_{sa}\times S)$, then there exists $W'' \in \op(X_{sa} \times S_{sa})$ with $W \towc W'' \towc W'$. $W$ is a finite union of $U_j \times V_j \in \op(X_{sa} \times S)$ with $U_j$ and $V_j$ relatively compact and $U_j \times \overline{V}_j \subset W'$. $\overline{V}_j$ being compact, it has a subanalytic neighborhood $V''_j$ such that $U_j \times V''_j \towc W'$. Then $W''=\cup_j (U_j \times V''_j) \in \op(X_{sa} \times S_{sa})$ and $W \towc W'' \towc W'$.

This implies that, given $W, W' \in \op(X_{sa} \times S)$ with $W \towc W'$ and $G \in \mod(\CC_{X_{sa} \times S})$, when we write
$$
\lind {W \towc W''} \Gamma(W'';G)
$$
we may assume that $W'' \in \op(X_{sa} \times S_{sa})$.

Let $F \in \mod(\CC_{X_{sa} \times S_{sa}})$ be $\T$-flabby. Then
$$
\lind {W \towc W''} \Gamma(W'';\imin aF) \simeq \lind {W \towc W''} \Gamma(W'';F)
$$
as a consequence of Proposition \ref{prop:imineta}. Since $F$ is $\T$-flabby the morphism $\Gamma(X \times S;F) \to \Gamma(W''; F)$ is surjective for each $W'' \in \op(X_{sa} \times S_{sa})$. So the morphism $\Gamma(X \times S;F) \to \lind {W \towc W''}\Gamma(W'';F) \simeq \lind {W \towc W''}\Gamma(W'';\imin aF)$ is surjective as well and factors through $\Gamma(W';\imin aF)$ by Proposition \ref{prop:imineta}. Then the result follows.
\end{proof}

\subsection{Construction of relative sheaves}

In the situation above, let  be given a sheaf $F$ on $(X\times S)_{sa}$.
As in \cite{MFP14}, we denote by  $F^{S,\sharp}$  the sheaf on $X_{sa} \times S_{sa}$ associated to the presheaf
\begin{eqnarray*}
\op(X_{sa} \times S_{sa})^{op} & \to & \mod(\CC) \\
U \times V & \mapsto & \Gamma(X \times V;\imin\rho\Gamma_{U \times S}F) \\
& \simeq & \Ho(\CC_U \boxtimes \rho_!\CC_V,F) \\
& \simeq & \lpro {W \subset\subset V, W\in\op(S_{sa})}\Gamma(U \times W;F).
\end{eqnarray*}

We set \begin{equation}\label{E:10}
F^{S}:=\imin a F^{S,\sharp}
\end{equation}
 and call it the relative sheaf associated to $F$. It is a sheaf on $X_{sa} \times S$.
 It is easy to check that $(\bullet)^S$ defines a left exact functor on $\mod(\CC_{(X \times S)_{sa}})$.
We will denote by $(\bullet)^{RS,\sharp}$ and $(\bullet)^{RS} \simeq \imin a \circ (\bullet)^{RS,\sharp}$ the associated right derived functors.

\begin{oss}\label{R1} When $X$ is reduced to a point, and $F$ is a sheaf on $S_{sa}$ then $F^S$ is nothing more than $\rho^{-1}F$ where $\rho$ is the morphism of sites $S\to S_{sa}$.

\end{oss}

\begin{lem} \label{lem:5}Let $U \in \op(X_{sa})$, $V \in \op(S)$ be relatively compact. Then
\begin{eqnarray*}
\Gamma(U \times V; F^{S}) & \simeq & \Gamma(X \times V;\imin\rho\Gamma_{U \times S}F) \\
& \simeq & \Ho(\CC_U \boxtimes \rho_!\CC_V,F).
\end{eqnarray*}
\end{lem}
\begin{proof} The second isomorphism follows by adjunction. Let us prove the first one. Let $U,V'$ be open subanalytic in $X$ and $S$ respectively. By Lemma 4.2 of \cite{MFP14} we have
$$\Gamma(U \times V';F^{S,\sharp}) \simeq \Gamma(X \times V';\imin\rho\Gamma_{U \times S}F) \simeq \lpro {W \subset\subset V'}\Gamma(U \times W;F).$$
By definition and Proposition \ref{prop:imineta}
$$\Gamma(U \times V;F^S) = \Gamma(U \times V;\imin aF^{S,\sharp}) \simeq \lpro {V' \subset\subset V}\Gamma(U \times V';F^{S,\sharp})$$
with $V'$ subanalytic. We have
$$
\lpro {V' \subset\subset V}\lpro {W \subset\subset V'}\Gamma(U \times W;F) \simeq \lpro {W \subset\subset V}\Gamma(U \times W;F) \simeq \Gamma(X \times V;\imin\rho\Gamma_{U \times S}F).
$$
This implies that $\Gamma(U \times V; F^{S}) \simeq \Gamma(X \times V;\imin\rho\Gamma_{U \times S}F)$ as required.
\end{proof}

\begin{prop} \label{prop:7} Let $G \in D^b(\CC_{X_{sa}})$ and $H \in D^b(\CC_S)$. Let $F \in D^b(\CC_{(X \times S)_{sa}})$.
Then
$$
\Rh(G \boxtimes H,F^{RS}) \simeq \Rh(G \boxtimes \rho_!H,F)$$ $$\simeq \Rh(\CC_X \boxtimes H,\imin\rho\rh(G \boxtimes \CC_S,F)) .
$$
\end{prop}
\begin{proof}
The second isomorphism follows by adjunction. We are going to prove the first one in several steps.

(a) Suppose that $F$ is injective. Then $F^{S,\sharp}$ is $\T$-flabby by Lemma 4.4 of \cite{MFP14} and $F^S=\imin a F^{S,\sharp}$ is c-soft by Proposition \ref{prop:Tflabbycsoft}. Let $U \in \op(X_{sa})$, let $V \in \op(S)$ and let us assume that both are relatively compact. We have
\begin{eqnarray*}
R\Gamma(U \times V;F^{RS}) & \simeq & \Gamma(U \times V;F^S) \\
& \simeq & \Ho(\CC_U \boxtimes \rho_!\CC_V,F) \\
& \simeq & \Rh(\CC_U \boxtimes \rho_!\CC_V,F),
\end{eqnarray*}
where the second isomorphism follows from Lemma \ref{lem:5}.

(b) Suppose that $G \in \mod^c_{\rc}(\CC_X)$ and $H \in \mod^c_{\rc}(\CC_S)$. Then $G$ (resp. $H$) is quasi-isomorphic to a bounded complex $G^\bullet$ (resp. $H^\bullet$) consisting of finite sums $\oplus \,\CC_W$ with $W$ subanalytic and relatively compact in $X$ (resp. $S$). Let $F \in D^b(\CC_{(X \times S)_{sa}})$ and let $F^\bullet$ be a complex of injective objects quasi-isomorphic to $F$. We have
\begin{eqnarray*}
\Rh(G \boxtimes H,F^{RS}) & \simeq & \Ho(G^\bullet \boxtimes H^\bullet,(F^\bullet)^S) \\
& \simeq & \Ho(G^\bullet \boxtimes \rho_!H^\bullet,F^\bullet) \\
& \simeq & \Rh(G \boxtimes \rho_!H,F),
\end{eqnarray*}
where the second isomorphism follows from (a).

(c) Suppose that $G \in \mod(\CC_{X_{sa}})$ and $H \in \mod^c_{\rc}(\CC_S)$. Then $G \simeq \lind i \,\rho_*G_i$, with $G_i \in \mod^c_{\rc}(\CC_X)$. We have
\begin{eqnarray*}
\Rh(\lind i\, \rho_*G_i \boxtimes H,F^{RS}) & \simeq & R\lpro i\,\Rh(G_i \boxtimes H,F^{RS}) \\
& \simeq & R\lpro i\,\Rh(G_i \boxtimes \rho_!H,F) \\
& \simeq & \Rh(\lind i\,\rho_*G \boxtimes \rho_!H,F),
\end{eqnarray*}
where for the first and the third isomorphism we refer to \cite{Pr99} and the second one follows from (b).

(d) Suppose that $G \in \mod(\CC_{X_{sa}})$ and $H \in \mod(\CC_S)$. Then $H \simeq \lind i\, H_i$, with $H_i \in \mod^c_{\rc}(\CC_X)$. We have
\begin{eqnarray*}
\Rh(G \boxtimes \lind i \,H_i,F^{RS}) & \simeq & R\lpro i\,\Rh(G \boxtimes H_i,F^{RS}) \\
& \simeq & R\lpro i\,\Rh(G_i \boxtimes \rho_!H_i,F) \\
& \simeq & \Rh(\lind i\,\rho_*G \boxtimes \rho_!\lind i\, H_i,F),
\end{eqnarray*}
where for the first and the third isomorphism we refer to \cite{Pr99} (we also used the fact that $\rho_!$ commutes with $\Lind$) and the second one follows from (c).

(e) Suppose that $G \in D^b(\CC_{X_{sa}})$ and $H \in D^b(\CC_S)$. Then $G$ (resp. $H$) is quasi-isomorphic to a bounded complex $G^\bullet$ (resp. $H^\bullet$) of sheaves on $X_{sa}$ (resp. $S$). By d\'evissage, we may reduce to the case $G \in \mod(\CC_{X_{sa}})$ and $H \in \mod(\CC_S)$ and the result follows from (d).

\end{proof}

\begin{cor} \label{cor:6 1/2} Let $G \in \mod(\CC_{X_{sa}})$, $H \in \mod(\CC_S)$ and let $F$ be an injective sheaf on $(X \times S)_{sa}$. Then $F^S$ is $\Ho(G \boxtimes H,\bullet)$-acyclic.
\end{cor}
\begin{proof}
It follows from Proposition \ref{prop:7} that $\Rh(G \boxtimes H,F^{RS}) \simeq \Rh(G \boxtimes \rho_!H,F)$ is concentrated in degree 0, as required.
\end{proof}

\begin{prop} \label{lem:9} Let $F \in D^b(\CC_{(X \times S)_{sa}})$. Let $G \in D^b(\CC_{X_{sa}})$ and $H \in D^b(\CC_S)$. Then
\begin{eqnarray*}
\imin {\rho_S} \rh(G \boxtimes H,F^{RS}) & \simeq  & \imin \rho \rh(G \boxtimes \rho_!H,F) \\
& \simeq & \rh(\CC_X \boxtimes H,\imin\rho\rh(G \boxtimes \CC_S,F)).
\end{eqnarray*}
In particular, when $G=\CC_X$ and $H=\CC_S$ we have $\imin\rho F \simeq \imin{\rho_S} F^{RS}\simeq \rho'^{-1} F^{RS,\sharp}$.
\end{prop}
\begin{proof} The second isomorphism follows by adjunction. Let us prove the first one.

(a) Let $K \in \mod(\CC_{(X \times S)_{sa}})$ and $K' \in \mod(\CC_{X_{sa} \times S})$. A morphism
$$
\eta_*K \to a_*H
$$
defines a morphism
$$
\imin a \eta_*K \to K'
$$
by adjunction and, composing with $\imin {\rho_S}$, a morphism
$$
\imin {\rho_S} \imin a \eta_*K \to \imin {\rho_S} K'.
$$
Remark that if $K \in \mod(\CC_{(X \times S)_{sa}})$ then $\imin\rho K \osi \imin {\rho_S} \imin a \eta_*K$. Indeed, for each $y \in X \times S$,
$$
(\imin\rho K)_y\simeq \lind {U \times V \ni y}K(U \times V) \simeq (\imin {\rho_S} \imin a \eta_*K)_y
$$
with $U \in \op(X_{sa})$, $V \in \op(S_{sa})$.

(b) Let us first suppose that $F,G,H$ are concentrated in degree zero. Hence, setting $K=\ho(G \boxtimes \rho_!H,F)$ and $K'= \ho(G \boxtimes H,F^S)$ in (a), to  any morphism $$ \eta_* \ho(G \boxtimes \rho_!H,F) \to a_* \ho(G \boxtimes H,F^S)$$ one associates a morphism $$\imin\rho \ho(G \boxtimes \rho_!H,F) \to \imin {\rho_S} \ho(G \boxtimes H,F^S).$$

Note that, for any $U, V\in\op((X\times S)_{sa})$, the natural morphism $\rho_!(H_V) \to (\rho_!H)_V$ induces a morphism $\Ho(G_U \boxtimes (\rho_!H)_V,F) \to \Ho(G_U \boxtimes \rho_!(H_V),F)$  hence a morphism $\psi: \eta_* \ho(G \boxtimes \rho_!H,F) \to a_* \ho(G \boxtimes H,F^S)$, which defines a morphism  $\imin\rho \ho(G \boxtimes \rho_!H,F) \to \imin {\rho_S} \ho(G \boxtimes H,F^S)$.

 Let us check on the  fibers that it is an isomorphism.
 Let $y \in X \times S$, then
\begin{eqnarray*}
(\imin\rho \ho(G \boxtimes \rho_!H,F))_y & \simeq & \lind {U \times V \ni y}\Ho(G_U \boxtimes (\rho_!H)_V,F) \\
& \simeq & \lind {U \times V \ni y}\,\lpro {W \subset\subset V}\Ho(G_U \boxtimes (\rho_!H)_W,F) \\
& \simeq & \lind {U \times V \ni y}\Ho(G_U \boxtimes \rho_!(H_V),F) \\
& \simeq & (\imin {\rho_S} \ho(G \boxtimes H,F^S))_y
\end{eqnarray*}
with $U \in \op(X_{sa})$, $V,W \in \op(S_{sa})$. The second and third isomorphisms follow since the morphism $(\rho_!H)_{V'} \to (\rho_!H)_{V}$ factors through $\ds{\lind {W \subset\subset V}(\rho_!H)_W}$ if $V'\subset\subset V$ ($V,V',W \in \op(S_{sa})$). Moreover, one has (see \cite{Pr08} for more details) that $\ds{\lind {W \subset\subset V}(\rho_!H)_W} \simeq \rho_!H \otimes \lind {W \subset\subset V}\rho_*\CC_W \simeq \rho_!H \otimes \rho_!\CC_V \simeq \rho_!(H_V)$.

(c) Suppose now that $F$ is injective and that $G,H$ are concentrated in degree 0. Let $U \in \op(X_{sa})$, $V \in \op(S_{sa})$. The complex
\begin{eqnarray*}
R\Gamma(U \times V;\rh(G \boxtimes H,F^S)) & \simeq & \Rh(G_U \boxtimes H_V,F^S)
\end{eqnarray*}
is concentrated in degree 0 by Corollary \ref{cor:6 1/2}. Then $F^S$ is $\shh om(G \boxtimes H,\bullet)$-acyclic.

(d) Let $G \in D^b(\CC_{X_{sa}})$ and $H \in D^b(\CC_S)$. Let $F \in D^b(\CC_{(X \times S)_{sa}})$ and let $F^\bullet$ be a complex of injective objects quasi-isomorphic to $F$. Then
\begin{eqnarray*}
\imin\rho \rh(G \boxtimes \rho_!H,F) & \simeq & \imin\rho \ho(G \boxtimes \rho_!H,F^\bullet) \\
& \simeq & \imin \rho_S \ho(G \boxtimes H,(F^\bullet)^S) \\
& \simeq & \imin \rho_S \rh(G \boxtimes H,F^{RS}),
\end{eqnarray*}
where the second isomorphism follows from (b) and the third one from (c).
\end{proof}

\begin{prop} \label{cor:8} Suppose that $F \in \mod(\rho_!\C^\infty_{X \times S})$ is $\Gamma(W;\bullet)$-acyclic for each $W \in \op((X \times S)_{sa})$. Then
$F$ is $(\bullet)^S$-acyclic.
\end{prop}
\begin{proof} Since $(\bullet)^{RS} \simeq \imin a \circ (\bullet)^{RS,\sharp}$, it is enough to show that $H^kF^{RS,\sharp}=0$ if $k \neq 0$. It is enough to prove that $F$ is $(\bullet)^{S,\sharp}$-acyclic. This follows from Corollary 4.9 of \cite{MFP14}.
\end{proof}

\begin{prop} \label{lem:8} Suppose that $F \in \mod(\rho_!\C^\infty_{X \times S})$ is $\Gamma(W;\bullet)$-acyclic for each $W \in \op((X \times S)_{sa})$. Then for each $U \in \op(X_{sa})$, $V \in \op(S)$ we have $R^k\Gamma(U \times V;F^{RS}) = 0$ if $k \neq 0$.
\end{prop}
\begin{proof} By Proposition \ref{lem:9}  we have $R\Gamma(U \times V;F^{RS}) \simeq R\Gamma(X \times V; \imin\rho R\Gamma_{U \times S}F)$. Since $F$ is $\Gamma(W;\bullet)$-acyclic for each $W \in \op((X \times S)_{sa})$, the complex $R\Gamma_{U \times S}F$ is concentrated in degree zero. Since $F$ is a $\rho_!\C^\infty_{X \times S}$-module,  $\imin\rho\Gamma_{U \times S}F$ is a $\C^\infty_{X \times S}$-module, hence c-soft and $\Gamma(X \times V;\bullet)$-acyclic. This shows the result. \end{proof}

\subsection{Behaviour of $(\bullet)^{RS}$ under pushforward and pull-back relatively to the parameter space}

We consider a real analytic map $h: S'\to S$ and still denote by $h$ the maps $\Id\times h: X\times S'\to X\times S$, $\Id\times h: X_{sa}\times S'\to X_{sa}\times S$.

\begin{prop}\label{change}
\begin{enumerate}
\item{Let $F\in D^b(\CC_{(X \times S)_{sa}})$.
Then there exists a natural morphism $$h^{-1}(F^{RS})\to (h^{-1}F)^{RS'}$$}
\item{

Let $G\in \rD^\rb(\CC_{(X\times S')_{sa}})$. Then there exists a natural isomorphism $$(Rh_*G)^{RS}\simeq Rh_*(G^{RS'})$$}
\end{enumerate}
\end{prop}

\begin{proof}

1) Let $U$ be a relatively compact subanalytic subset of $X$ and let $V'$ be any open set in $S'$. Let $V$ be any open subset of $S$ such that $V\supset h(V')$. By Proposition \ref{lem:9}  we have

\begin{eqnarray*}
R\Gamma(U \times V; F^{RS}) & \simeq & R\Gamma(X \times V; \imin\rho R\Gamma_{U \times S}F) \\
& \to & R\Gamma(X \times V'; \imin\rho h^{-1}R\Gamma_{U \times S}F) \\
& \to & R\Gamma(X \times V';\rho^{-1}R\Gamma_{U\times S'}(h^{-1}F)) \\
& \simeq & R\Gamma(U \times V'; (h^{-1}F)^{RS'})
\end{eqnarray*}
 which gives the desired morphism.

2) Given $U$ and $V$ as in 1), let us set $V'=h^{-1}(V)$. Then $$R\Gamma (U\times V; (Rh_*G)^{RS})\simeq R\Gamma (X\times V; \rho^{-1}R\Gamma_{U\times S}Rh_*G)$$$$\simeq R\Gamma (X\times V; \rho^{-1}Rh_*R\Gamma_{U\times S'}G)\simeq R\Gamma (X\times V; Rh_*\rho^{-1}R\Gamma_{U\times S'}G)$$$$\simeq R\Gamma(X\times V'; \rho^{-1}R\Gamma_{U\times S'}G)\simeq R\Gamma (U\times V'; G^{RS'})$$ $$\simeq R\Gamma (U\times V; Rh_*G^{RS'})$$
where all the isomorphisms commute with restrictions. This gives the desired isomorphism.
\end{proof}

\subsection{The S-locally constant case}

The result below is an improvement of Proposition 3.3 of \cite{MFCS2}. The latter was a step in the construction of the relative Riemann-Hilbert functor, denoted by $\RH_X^S$ where $d_S=1$ was assumed. As already said, there the construction relied on the possibility of considering on $S_{sa}$ a basis formed by Stein open subanalytic sets. This difficulty is overcome with the new site $X_{sa}\times S$ and the proof goes then in a similar way.

\begin{lem} \label{lem:O_S acyclic} The sheaf $\imin p\OO_S$ is $\rho_{S*}$-acyclic.
\end{lem}
\begin{proof} Suppose that $U \in \op(X_{sa})$ is connected and $V \in \op(S)$ is Stein. Then, by Proposition 3.3.9 of \cite{KS1}, we have $$R\Gamma(U \times V;R\rho_{S*}\imin p \OO_S) \simeq R\Gamma(U \times V;\imin p \OO_S) \simeq R\Gamma(V;\OO_S)$$  and the latter is concentrated in degree 0. The result follows since the topology of $X_{sa} \times S$ admits a basis consisting of open subsets $U \times V$ with $U \in \op(X_{sa})$ contractible (as a consequence of \cite{Wi05}), and $V$ Stein in $S$.
\end{proof}

\begin{prop}\label{locconst2} If $F \in \mod(\CC_{X \times S})$ is locally isomorphic to $p^{-1}G$, where $G$ is a coherent  $\sho_S$-module, then $R\rho_{S*}F\simeq \rho_{S!}F$.
\end{prop}

\begin{proof} We have to show that, if $F \in \mod(\CC_{X \times S})$ is locally isomorphic to $p^{-1}G$, where $G$ is $\OO_S$-coherent, the natural morphism
$$
\rho_{S!}F \to R\rho_{S*}F,
$$
defined by the adjunction morphism $\id \to \imin{\rho_S}R\rho_{S*} \simeq \id$ is an isomorphism. It is enough to check it on a basis for the topology of
$X_{sa} \times S$. First of all, let us prove that $R\rho_{S*}F$ is concentrated in degree zero. There exists a covering $\mathcal{U}$ of $X_{sa} \times S$ such that, for each $W \in \mathcal{U}$, one has $(R\rho_{S*}F)|_W \simeq   (R\rho_{S*}p^{-1}G)|_W$, where $G$ is $\OO_S$-coherent. Then Lemma \ref{lem:O_S acyclic} implies that $R\rho_*F$ is concentrated in degree zero.

Up to taking a refinement, we may suppose that $W=U \times V$ with $U \in \op(X_{sa})$ relatively compact and connected and $V \in \op(S)$ connected. In this case $\rho_{S!}p^{-1}G$ is the sheaf associated to the presheaf $$U \times V \mapsto \lind {U'}\,\Gamma(U' \times V;p^{-1}G)$$ with $ {\overline{U} \subset U'}$. Moreover, if $U$ is connected, we can suppose that $U'$ is connected as well. Indeed, if $U$ is subanalytic, every subanalytic neighborhood of $\overline{U}$ contains a subanalytic neighborhood $U'$ which is connected (as an easy consequence of the triangulation of subanalytic sets).

The result follows since, according to Proposition 3.3.9 of \cite{KS1}, $\Gamma(U \times V; p^{-1}G) \simeq \Gamma(V;G)$.
\end{proof}

\subsection{The  sheaves  $\C_{X\times S}^{\infty, t, S}$, $\db_{X\times S}^{t,S}$, $\C_{X\times S}^{\infty,\mathrm{w},S}$, $\OO_{X\times S}^{t,S}$ and $\OO_{X\times S}^{\mathrm{w},S}$}\label{S:3}
Let $X$  and $S$ be real analytic manifolds.
The construction given by (\ref{E:10}) allows us to introduce the following sheaves:
\begin{enumerate}
\item{$\C^{\infty,t,S}_{X\times S}:=(\C^{\infty,t}_{X\times S})^S$ as the relative sheaf associated to $\C^{\infty,t}_{X\times S}$,}
\item{$\db_{X\times S}^{t,S}:=(\dbt_{X \times S})^S$ as the relative sheaf associated  to $\db^{t}_{X\times S}$,}
\item{$\C_{X\times S}^{\infty,\mathrm{w},S}:=(\CW_{X \times S})^S$ as the relative sheaf associated to $\C_{X\times S}^{\infty,\mathrm{w}}$.}
\end{enumerate}

Let us apply the results of the previous section to these sheaves.

\begin{prop} \label{prop:6}Let $U \in \op(X_{sa})$, $V \in \op(S_{sa})$. Then
\begin{enumerate}
\item{$\Gamma(U \times V;\C^{\infty, t,S}_{X\times S})\simeq\Gamma(X \times V;\imin\rho\Gamma_{U \times S}\C^{\infty,t}_{X \times S}) \\ \ \ \ \ \simeq \Gamma(X \times V;T\ho(\CC_{U \times S},\C^\infty_{X \times S})),$}
\item{$\Gamma(U \times V;\db^{t,S}_{X \times S})\simeq\Gamma(X \times V;\imin\rho\Gamma_{U \times S}\dbt_{X \times S}) \\ \ \ \ \ \simeq \Gamma(X \times V;T\ho(\CC_{U \times S},\db_{X \times S})),$}
\item{$\Gamma(U \times V;\C^{\infty,\mathrm{w},S}_{X \times S})\simeq\Gamma(X \times V;\imin\rho\Gamma_{U \times S}\C^{\infty,\mathrm{w}}_{X \times S}) \\ \ \ \ \ \simeq \Gamma(X \times V;H^0D'\CC_U \boxtimes \CC_S \wtens \C^\infty_{X \times S}).$}
\end{enumerate}
\end{prop}

\begin{oss} In the case of subanalytic open sets like $U \times V$ the definition of
the sheaf $\db^{t,S}_{X \times S}$ is well known. A section $s \in \Gamma(U \times V ;\db_{X \times S})$ belongs
to $\Gamma(U \times V ;\db^{t,S}_{X \times S})$ if it extends as a section of $\Gamma(X \times V ;\db_{X \times S})$. It is a consequence of the definition of the functor $T\ho$ of \cite{Ka3}.
\end{oss}

\begin{prop} \label{prop:8}

i) Suppose that $\mathcal{F} =\db^t_{X \times S},\C^{\infty,t}_{X \times S}, \C^{\infty,\mathrm{w}}_{X \times S}$. Then $\mathcal{F}$ is $(\bullet)^S$-acyclic. Moreover $\db^{t,S}_{X \times S},\C^{\infty,t,S}_{X \times S}$ are $\Gamma(U \times V;\bullet)$-acyclic for each $U \in \op(X_{sa})$, $V \in \op(S)$.

ii)  $\C^{\infty,\mathrm{w},S}_{X \times S}$ is $\Gamma(U \times V;\bullet)$-acyclic for each $U \in \op(X_{sa})$ locally cohomologically trivial and $V \in \op(S)$.
\end{prop}

\begin{prop} \label{prop 7}Let $G \in D^b(\CC_{X_{sa}})$, $H \in D^b(\CC_S)$. Then
\begin{enumerate}
\item{$\imin\rho_S\rh(G \boxtimes H,\C^{\infty,t,S}_{X \times S}) \simeq \imin\rho\rh(G \boxtimes \rho_!H,\C^{\infty,t}_{X \times S}) \\ \ \ \ \ \simeq \rh(\CC_X \boxtimes H,\imin\rho\ho(G \boxtimes \CC_S,\C^{\infty,t}_{X \times S})),$}
\item{$\imin\rho_S\rh(G \boxtimes H,\db^{t,S}_{X \times S}) \simeq \imin\rho\rh(G \boxtimes \rho_!H,\dbt_{X \times S}) \\ \ \ \ \ \simeq \rh(\CC_X \boxtimes H,\imin\rho\ho(G \boxtimes \CC_S,\dbt_{X \times S})),$}
\item{$\imin\rho_S\rh(G \boxtimes H,\C^{\infty,\mathrm{w},S}_{X \times S}) \simeq \imin\rho\rh(G \boxtimes \rho_!H,\CW_{X \times S}) \\ \ \ \ \ \simeq \rh(\CC_X \boxtimes H,\imin\rho\ho(G \boxtimes \CC_S,\CW_{X \times S})).$}
\end{enumerate}
When $G \in D^b_{\rc}(\CC_X)$ we have
\begin{enumerate}
\item{$\imin\rho_S\rh(G \boxtimes H,\C^{\infty,t,S}_{X \times S})  \simeq \rh(\CC_X \boxtimes H,T\ho(G \boxtimes \CC_S,\C^\infty_{X \times S})),$}
\item{$\imin\rho_S\rh(G \boxtimes H,\db^{t,S}_{X \times S}) \simeq \rh(\CC_X \boxtimes H,T\ho(G \boxtimes \CC_S,\db_{X \times S})),$}
\item{$\imin\rho_S\rh(G \boxtimes H,\C^{\infty,\mathrm{w},S}_{X \times S}) \simeq \rh(\CC_X \boxtimes H,D'G \boxtimes \CC_S \wtens \C^\infty_{X \times S}).$}
\end{enumerate}
In particular, when $G=\CC_X$ and $H=\CC_S$ we have $\imin\rho_S\C^{\infty,t,S}_{X \times S} \simeq \C^\infty_{X \times S}$, $\imin\rho_S\db^{t,S}_{X \times S} \simeq \db_{X \times S}$, $\imin\rho_S\C^{\infty,\mathrm{w},S}_{X \times S} \simeq \C^\infty_{X \times S}$.
\end{prop}

Now we are going to study the action of differential operators on these new sheaves.

\begin{lem}\label{L:3}
There is a natural action of $\rho_{S!}\D_{X\times S}$ on $\C^{\infty,t,S}_{X\times S}$ and on $\C_{X\times S}^{\infty,\mathrm{w},S}$.
\end{lem}
\begin{proof}
Let $K \in \mod(\CC_{(X \times S)_{sa}})$. In order to prove the action of $\rho_{S!}\D_{X\times S}$ on $K^S$ it is enough to prove the action of $\rho'_!\D_{X \times S}$ on $K^{S,\sharp}$. When $K=\db^{t,S}_{X\times S},\C^{\infty,t,S}_{X\times S},\C_{X\times S}^{\infty,\mathrm{w},S}$, it is a consequence of Lemma 5.4 of \cite{MFP14}.
\end{proof}

Let us now assume that $X$ and $S$  are complex manifolds and consider the projection $p:X\times S\to S$. Let us denote as usual by $\overline{X}\times\overline{S}$ the complex conjugate manifold. Identifying the underlying real analytic manifold $X_{\R} \times S_{\R}$ to the diagonal of $(X \times S) \times (\overline{X} \times \overline{S})$, we have:

\begin{lem}\label{lem:O_S action}
 $\rho_{S*}\imin p\OO_S$ acts on $\db^{t,S}_{X \times S}$, on $\C^{\infty,t,S}_{X \times S}$ and on $\C_{X \times S}^{\mathrm{w},S}$.
\end{lem}
\begin{proof}
Let $K=\db^{t}_{X\times S},\C^{\infty,t}_{X\times S},\C_{X\times S}^{\infty,\mathrm{w}}$. By Lemma 5.5 of \cite{MFP14}, $\rho'_*\imin p \OO_S$ acts on $K^{S,\sharp}$. Note that $\rho'_* \simeq a_* \circ \rho_{S*}$, then $\imin a \rho'_* \imin p \OO_S \simeq \rho_{S*}\imin p \OO_S$ acts on $\imin a K^{S,\sharp} = K^S$ as required.
\end{proof}

We also have natural actions of
$\rho_{S!}\D_{\overline{X}\times \overline{S}}, \rho_{S!}\D_{\overline{X}\times S}$ on $\db^{t,S}_{X\times S}$.

The construction given by (\ref{E:10}) allows us to introduce the following objects of $D^b(\CC_{X_{sa} \times S})$:
\begin{enumerate}
\item{$\OO_{X \times S}^{t,S}:=(\ot_{X \times S})^{RS}$, the relative sheaf associated  to $\ot_{X \times S}$, that is $$\OO_{X \times S}^{t,S}  \simeq  (\rh_{\rho_!\D_{\overline{X} \times \overline{S}}}(\rho_!\OO_{\overline{X} \times \overline{S}},\dbt_{X \times S}))^{RS}$$ $$ \simeq  (\rh_{\rho_!\D_{\overline{X} \times \overline{S}}}(\rho_!\OO_{\overline{X} \times \overline{S}},\C^{\infty,t}_{X \times S}))^{RS}$$}
\item{$\OO_{X \times S}^{\mathrm{w},S}:=(\OW_{X \times S})^{RS}$, the relative sheaf associated to $\OW_{X \times S}$, that is $$\OO_{X \times S}^{\mathrm{w},S} \simeq  (\rh_{\rho_!\D_{\overline{X} \times \overline{S}}}(\rho_!\OO_{\overline{X} \times \overline{S}},\CW_{X \times S}))^{RS}.$$}
\end{enumerate}
The  exactness of $\rho_{S!}$   together with Proposition
\ref{prop:8}
 allow to conclude:
\begin{prop} We have the following isomorphisms in $D^b(\CC_{X_{sa} \times S})$.
%\begin{eqnarray*}
$$\OO_{X \times S}^{t,S} \simeq \rh_{\rho_{S!}\D_{\overline{X} \times \overline{S}}}(\rho_{S!}\OO_{\overline{X} \times \overline{S}},\db^{t,S}_{X \times S})$$ $$
 \simeq \rh_{\rho_{S!}\D_{\overline{X} \times \overline{S}}}(\rho_{S!}\OO_{\overline{X} \times \overline{S}},\C^{\infty,t,S}_{X \times S})$$
$$\OO_{X \times S}^{\mathrm{w},S} \simeq \rh_{\rho_{S!}\D_{\overline{X} \times \overline{S}}}(\rho_{S!}\OO_{\overline{X} \times \overline{S}},\C^{\infty,\mathrm{w},S}_{X \times S}).$$
%\end{eqnarray*}
\end{prop}

Proposition \ref{lem:9} together with Proposition 7.3.2 of \cite{KS5} entail:
\begin{prop} \label{prop:7BIS}Let $G \in D^b(\CC_{X_{sa}})$, $H \in D^b(\CC_S)$. Then
\begin{enumerate}
\item{$\imin\rho_S\rh(G \boxtimes H,\OO^{t,S}_{X \times S}) \simeq \imin\rho\rh(G \boxtimes \rho_{!}H,\ot_{X \times S}) \\ \ \ \ \ \simeq \rh(\CC_X \boxtimes H,\imin\rho\rh(G \boxtimes \CC_S,\ot_{X \times S})),$} \\
\item{$\imin\rho_S\rh(G \boxtimes H,\OO^{\mathrm{w},S}_{X \times S}) \simeq \imin\rho\rh(G \boxtimes \rho_{!}H,\OW_{X \times S}) \\ \ \ \ \ \simeq \rh(\CC_X \boxtimes H,\imin\rho\rh(G \boxtimes \CC_S,\OW_{X \times S})).$}
\end{enumerate}
When $G \in D^b_{\rc}(\CC_X)$ we have
\begin{enumerate}
\item{$\imin\rho_S\rh(G \boxtimes H,\sho^{t,S}_{X \times S})
 \simeq \rh(\CC_X \boxtimes H,T\ho(G \boxtimes \CC_S,\sho_{X \times S})),$} \\
\item{$\imin\rho_S\rh(G \boxtimes H,\sho^{\mathrm{w},S}_{X \times S})
 \simeq \rh(\CC_X \boxtimes H,D'G \boxtimes \CC_S \wtens \sho_{X \times S}).$}
\end{enumerate}
In particular, when $G=\CC_X$ and $H=\CC_S$ we have $\imin\rho_S\sho^{t,S}_{X \times S} \simeq \sho_{X \times S}$, $\imin\rho_S\sho^{\mathrm{w},S}_{X \times S} \simeq \sho_{X \times S}$.
\end{prop}

The examples given in \cite{MFP14} can now be stated in a more general case ($V$ needs no longer to be subanalytic):
\begin{es}
Let $U=\{z\in\CC, \Im z>0\}$, let $V$ be open in $\CC^n$ and let $g(s)$ be a holomorphic function on $V$. Then, after a choice of a determination of $log\,z$ on $U$, $z^{g(s)}$ defines a section of $\Gamma(U\times V; \OO^{t,S}_{\CC\times \CC^n})$.
\end{es}

\begin{es}Let $U=\R_{>0}$ with a coordinate $x$, let $V$ be an open set in $\R$ and let $a(s)$ be any continuous function on $V$. Let $f\in\Gamma(\Omega\setminus V;\OO_{\CC})$, where $\Omega$ is an open neighborhood of $V$ in $\CC$, be such that $a$ is the boundary value $vb(f)$ of $f$ as a hyperfunction. Then $x^{a}_{+}:=vb(z^f)$, with $arg \,z\in ]0,2\pi[$,  is a section of $\Gamma(U\times V; \db^{t,S}_{\R\times \R})$.
\end{es}

\subsection{Functorial properties on the parameter space}

\begin{oss}\label{R2} We remark that, for any $S$, the site $X_{sa}\times S$ is a ringed site both relatively to the sheaf $\rho_{S*}(p^{-1}\sho_S)$ and to the sheaf $\rho_{S!} \sho_{X\times S}$ (cf \cite{KS3}, page 449).
Given a morphism of complex manifolds $\pi: S'\to S$,
we have isomorphisms of functors $$\rho_{S'}^{-1}\pi^{-1}\rho_{S*}\simeq \pi^{-1}$$ $$\rho_{S'*}\rho_{S'}^{-1}\pi^{-1}\rho_{S*}\simeq \rho_{S'*}\pi^{-1}$$ and composing with the natural morphism $\Id\to \rho_{S'*}\rho_{S'*}^{-1}$ gives a natural morphims
$$\pi^{-1}\rho_{S*}\to \rho_{S'*}\pi^{-1}$$

Since $\sho_{S'}$ is a $\pi^{-1}\sho_S$-module, $\rho_{S'*}(p^{-1}\sho_S')$ is a $\rho_{S'*}(p^{-1}\pi^{-1}\sho_S)$-module hence a $\pi^{-1}\rho_{S*}(p^{-1}\sho_S)$-module. Similarly, $\rho_{S'!}\sho_{X\times S'}$ is a $\pi^{-1}\rho_{S!}\sho_{X\times S}$-module.

In other words, $\pi$ induces a morphism of ringed sites with respect to both sheaves of rings.
%$\sho_{X\times S}^{t,S}$ is a $\rho_{S*}\sho_{X\times S}$-module.

Consequently, according to loc.cit Theorem 18.6.9 (i),  the derived functors $L\pi^*: D(\rho_{S*}p^{-1}\sho_S)\to D(\rho_{S'*}p^{-1}\sho_{S'})$ resp. (keeping the same notation $\pi$ for the morphism $\Id\times \pi$), $L\pi^*: D(\rho_{S!}\sho_{X\times S})\to D(\rho_{S'!}\sho_{X\times S'})$ are well defined.
\end{oss}

If $S'$ is a closed submanifold of $S$, we denote by $i_{S'}$ the closed immersion
$i_{S'}: X\times S'\hookrightarrow X\times
S.$  Hence $Li^*_{S'}$ is the functor on $D^{b}(\rho_{S\ast}(p^{-1}\sho_S))$ given by $\shf\mapsto i^{-1}_{S'}(\rho_{S\ast}(p^{-1}(\sho_S/\mathfrak{m}))\otimes^L_{\rho_{S\ast}(p^{-1}\sho_S)}\shf)$ where $\mathfrak{m}$ is the sheaf of ideals of $\sho_S$ of functions vanishing on $S'$.

We denote by $\bD'$ the functor on $D^b(p^{-1}\sho_S)$ given by
$$\bD'(F)=\rh_{p^{-1}\sho_S}(F, p^{-1}\sho_S)$$

\begin{prop}[Action by a closed immersion]\label{P:Lis}
Let $S'$ be a closed submanifold of $S$. There are natural morphisms in $D^{b}(\rho_{S\ast}(p^{-1}\sho_S))$
\begin{enumerate}
\item{$Li^*_{S'}\sho^{t,S}_{X \times S}\to \sho^{t,S'}_{X\times S'}$}
\item{$Li^*_{S'}\OO_{X \times S}^{\mathrm{w},S} \to \OO_{X\times S' }^{\mathrm{w,S'}}$}

\end{enumerate}
which are isomorphisms.
\end{prop}

\begin{proof}
Let us construct the morphism $(1)$. (Note that for $\codim S'=0$ the result is trivial.)

This amounts to showing that the following morphism is an isomorphism
$$\Gamma(V; \sho_S/\mathfrak{m})\otimes^L_{\Gamma(V; \sho_S)}R\Gamma(U\times V; \sho^{t,S}_{X\times S})\to R\Gamma(U\times W; \sho^{t,S'}_{X\times S'})$$ for any relatively compact open subset $U\subset X$ and any open subset $V\subset S$ running on a basis of the topology of $S$ consisting of Stein open relatively compact subsets and where we note $W=V\cap S'$. We note that, according to $(2)$ of Proposition \ref{prop:6}, we have
$$R\Gamma(U\times V; \sho^{t,S}_{X\times S})=R\Gamma(X\times V; T\ho(\CC_{U\times S}, \sho_{X\times S}))$$

 Recall that, for a $\shd_{X\times S}$-module $\shm$, the pull-back by $i_{S'}$ is given by $_Di_{S'}^*\shm:=\sho_{X\times S'}\otimes^L_{i^{-1}\sho_{X\times S}}i_{S'}^{-1}\shm$. Note that $\sho_{X\times S}$ is a flat $p^{-1}\sho_S$-module. Hence $$\sho_{X\times S'}=i_{S'}^{-1}(\sho_{X\times S}/\mathfrak{m}\sho_{X\times S})$$ $$=i_{S'}^{-1}(p^{-1}(\sho_S/\mathfrak{m})\otimes_{p^{-1}\sho_S}\sho_{X\times S})$$

Therefore, Theorem 5.8 (5.16) of \cite{KS4}  applied to the closed embedding $i_{S'}$ and to $F=\CC_{U\times S}$ becomes a natural isomorphism in $D^b(\shd_{X\times S'})$
\begin{equation}
\label{eq:$(A)$} T\ho(\CC_{U\times S'},\sho_{X\times S'})\simeq Li^*_{S'}T\ho(\CC_{U\times S}, \sho_{X\times S}).
\end{equation}
Up to a shrinking of $V$ so that a family of holomorphic coordinates vanishing on $S'$ is defined in $V$, by taking the corresponding Koszul resolution ($\sho_S$-free) of $\sho_S/\mathfrak{m}$, we have

$$\Gamma(V; \sho_S/\mathfrak{m})\otimes^L_{\Gamma(V; \sho_S)}R\Gamma(X\times V; T\ho(\CC_{U\times S}, \sho_{X\times S}))$$ $$\simeq R\Gamma(X\times V; p^{-1}(\sho_S/\mathfrak{m})\otimes^L_{p^{-1}\sho_S}T\ho(\CC_{U\times S}, \sho_{X\times S}))$$ and this last term, according to \eqref{eq:$(A)$}, is isomorphic to

$
R\Gamma(X\times W; T\ho(\CC_{U\times S'}, \sho_{X\times S'}))
\simeq
R\Gamma(U\times W; \sho^{t,S'}_{X\times S'})
$

\noindent hence the desired result follows since these isomorphisms are compatible with restrictions to open subsets.

The construction of morphism $(2)$ is similar using ($3$) of Proposition \ref{prop:6} and (5.15) of Theorem 5.8  of \cite{KS4}.
\end{proof}

\begin{prop}\label{rhomot}
Let $S'\subset S$ be a closed submanifold of $S$ of codimension $d$. Then we have an isomorphism in $D^b(\shd_{X\times S'/S'})$
$$i_{S'}^{-1}\rh_{\rho_{S*}p^{-1}\sho_S}(\rho_{S*}p^{-1}(\sho_S/\mathfrak{m}), \sho^{t,S}_{X\times S})
\simeq \sho^{t,S'}_{X\times S'}[-d]$$
\end{prop}
\begin{proof}
We have $$\rh_{\rho_{S*}p^{-1}\sho_S}(\rho_{S*}p^{-1}(\sho_S/\mathfrak{m}), \sho^{t,S}_{X\times S})
\simeq \rho_{S*}\bD'(p^{-1}(\sho_S/\mathfrak{m}))\otimes^L_{\rho_{S*}p^{-1}\sho_S}\sho^{t,S}_{X\times S}$$  $$\simeq \rho_{S*}(p^{-1}(\sho_S/\mathfrak{m}))[-d]\otimes^L_{\rho_{S*}p^{-1}\sho_S} \sho^{t,S}_{X\times S}$$ where the first isomorphism is given by Lemma 3.22 of \cite{MFCS2}.

Applying $i_{S'}^{-1}$ we derive un isomorphism with $
Li^*_{S'}\sho^{t,S}_{X\times S}[-d]$.
The result then follows by Proposition \ref{P:Lis}.
\end{proof}

As a consequence of Proposition \ref{P:Lis} we conclude:
\begin{prop}[Restriction to the fibers]\label{P:Lis2}
Let $s_0\in S$. We simply denote by $Li^*_{s_0}$ the functor $Li^*_{\{s_0\}}$ on $D^{b}(\rho_{S\ast}(p^{-1}\sho_S))$.
There are natural morphisms
\begin{enumerate}
\item{$Li^*_{s_0}\sho^{t,S}_{X \times S}\to \sho^t_X$}
\item{$Li^*_{s_0}\OO_{X \times S}^{\mathrm{w},S} \to \OO_{X }^{\mathrm{w}}$}
\end{enumerate}
which are isomorphisms, where we  identify $X_{\sa}$ with $X_{\sa}\times \{s_0\}$.
\end{prop}

More generally we consider now a morphism $\pi: S'\to S$ of complex analytic manifolds.
\begin{prop}[Inverse image]\label{Pinv}
There exists a natural morphism
$$ \pi^{-1}\sho^{t,S}_{X\times S}\to\sho^{t,S'}_{X\times S'}$$

which coincides with the usual morphism

 $$\pi^{-1}\sho_{S}\to\sho_{S'}$$ when $X=pt$. In particular we have a natural morphism $L\pi^*\sho^{t.S}_{X\times S}\to \sho^{t,S'}_{X\times S'}$.

\end{prop}

\begin{proof}

We recall that the case of a ramification of finite degree was already treated in Lemma 2.11 of \cite{FMFS} in a different but equivalent framework (as said above, there one considered the relative site $X_{sa}\times S_{sa}$ instead of $X_{sa}\times S$).
The same argument allows us to reduce to prove the existence of a morphism in $D^b(\rho_{S'*}\pi^{-1}\shd_{X\times S/S})$ $$\pi^{-1}C_{X\times S}^{\infty, t, S}\to C_{X\times S'}^{\infty,t,S'}$$
In that case the morphism is nothing more than the composition with $\pi$ which commutes with operators in $\pi^{-1}\DXS$ and keeps the growth conditions.

\end{proof}

\subsection{Application}

Let $Y$ be an hypersurface of the manifold $X$ and assume that $Y$ is a normal crossing divisor. We shall note $X^*:=X\setminus Y$ and denote by $j$ either the open inclusion $j:X^*\to X$ or $j: X^*\times S\to X\times S$. We also keep $j$ to denote the associated morphism $j: X_{sa}^*\times S\to X_{sa}\times S$. We still denote by $p$ the restriction of $p$ to $X^*\times S$.

Recall that the notion of $S$-local system on $X^*\times S$ goes back to the work of Deligne (\cite{De}) and was the object of a systematic study in the Appendix of \cite{MFCS2}.

A sheaf of  $p^{-1}\sho_S$-modules is $S$-locally constant (or an $S$-local system for short) if there exists a coherent $\sho_S$-module $G$ such that, locally on $X^*\times S$, $F\simeq p^{-1}G$.

The following result is an improvement of a statement in Lemma 3.25 of \cite{MFCS2} where the locally free case was considered assuming $d_S=1$.
 It is one important step to the relative Riemann-Hilbert correspondence for arbitrary $d_S$.
\begin{prop}\label{RHS}
Let $F$ be an $S$-local system on $X^*\times S$. Then $$\shm:=\rho_S^{-1}Rj_*(\rho_{S^*}F\otimes^L_{\rho_{S*}p^{-1}\sho_S}j^{-1}\sho_{X\times S}^{t,S})$$ is concentrated in degree zero and $$\shh^0 \shm\simeq \rho_S^{-1}j_*(\rho_{S*}F\otimes_{\rho_{S*}p^{-1}\sho_S}j^{-1}\sho_{X\times S}^{t,S})$$

\end{prop}

\begin{proof}

Firstly we remark that the statements are local in $X\times S$.

We start by proving the first statement.
Let $G$ be a coherent $\sho_S$-module such that $F$ is locally isomorphic to $p^{-1}G$.
We have
$$Rj_*(\rho_{S^*}F\otimes^L_{\rho_{S*}p^{-1}\sho_S}j^{-1}\sho_{X\times S}^{t,S})$$ $$\simeq
Rj_*(\rho_{S*}F\otimes^L_{\rho_{S*}p^{-1}\sho_S}(\rho_{S_!}\sho_{X^*\times S}\otimes^L_{\rho_{S!}\sho_{X^*\times S}}j^{-1}\sho_{X\times S}^{t,S}))$$ $$\simeq Rj_*((\rho_{S_!}F\otimes^L_{\rho_{S_!}p^{-1}\sho_S}\rho_{S_!}\sho_{X^*\times S})\otimes^L_{\rho_{S!}\sho_{X^*\times S}}j^{-1}\sho_{X\times S}^{t,S})$$ thanks to Proposition \ref{locconst2}.

Since $\sho_{X^*\times S}$ is $p^{-1}\sho_S$-flat and $\rho_{S_!}$ is exact and commutes with $\otimes$,  we have $$\rho_{S_!}F\otimes^L_{\rho_{S_!}p^{-1}\sho_S}\rho_{S_!}\sho_{X^*\times S}\simeq \rho_{S_!}(F\otimes_{p^{-1}\sho_S}\sho_{X^*\times S})$$

According to the proof of Theorem 2.6 of \cite{MFCS2} we have an isomorphism of $\sho_{X^*\times S}$-modules
$$F\otimes_{p^{-1}\sho_S}\sho_{X^*\times S}\simeq p^{-1}G\otimes_{p^{-1}\sho_S}\sho_{X^*\times S}$$ hence we may assume from the beginning that $F=p^{-1}G$.

Therefore we conclude
$$\shm\simeq \rho_S^{-1}Rj_*(\rho_{S*}j^{-1}p^{-1}G\otimes_{\rho_{S*}p^{-1}\sho_S} j^{-1}\sho_{X\times S}^{t,S})$$
$$\simeq \rho_S^{-1}Rj_*(j^{-1}\rho_{S*}p^{-1}G\otimes_{\rho_{S*}p^{-1}\sho_S} j^{-1}\sho_{X\times S}^{t,S})$$
$$\simeq \rho_S^{-1}(\rho_{S*}p^{-1}G\otimes_{\rho_{S*}p^{-1}\sho_S}Rj_* j^{-1}\sho_{X\times S}^{t,S})$$
$$\simeq p^{-1}G\otimes_{p^{-1}\sho_S}\rho_S^{-1}Rj_* j^{-1}\sho_{X\times S}^{t,S}$$ where, in the last isomorphism, we used the commutation of $\rho_S^{-1}$ with $\otimes$.

The statement follows because, according to Proposition \ref{prop:7}, we have  $\rho_S^{-1}Rj_* j^{-1}\sho_{X\times S}^{t,S}\simeq T\ho(\CC_{X^*\times S}, \sho_{X\times S})\simeq\sho_{X\times S}(*Y\times S)$

It remains to prove that, for a basis $(U\times V)$ of the topology of $X\times S$, which can be chosen so that $U$ and $V$ are Stein open sets and $U$ subanalytic relatively compact,  $R\Gamma(U\times V;  Rj_* j^{-1}\sho_{X\times S}^{t,S})$ is concentrated in degree zero.

We have $$R\Gamma(U\times V;  Rj_* j^{-1}\sho_{X\times S}^{t,S})\simeq R\Gamma((U\setminus Y)\times V; \sho_{X\times S}^{t,S})\simeq R\Gamma((U\setminus Y)\times V; \sho^t_{X\times V})$$ and by Corollary 2.3.5 of \cite{S21} this last complex is concentrated in degree zero and thus equals $\Gamma(U\setminus Y)\times V; \sho^t_{X\times V})$ which ends the proof.
\end{proof}

\author{Teresa Monteiro Fernandes}\\
{Centro de Matem\'atica, Aplica\c{c}\~{o}es Funda\-men\-tais e Investiga\c c\~ao operacional and Departamento de Matem\' atica da Faculdade de Ci\^en\-cias da Universidade de Lisboa, Bloco C6, Piso 2, Campo Grande, 1749-016, Lisboa,
Portugal}\\
{mtfernandes@fc.ul.pt}

\author{Luca Prelli}\\
{Dipartimento di Matematica ``Tullio Levi-Civita'' Universit\`a degli Studi di Padova,
Via Trieste, 63,
35121 Padova Italy}\\
{lprelli@math.unipd.it}

\end{document}